\documentclass{elsart}

\usepackage{natbib}

\usepackage{amssymb}

\usepackage[french,english]{babel}

\usepackage{amsmath,amscd,amsfonts,euscript}
\usepackage{enumerate}
\usepackage{theorem}
\usepackage{textcomp}

\newtheorem{theo}{Theorem}[section]
\newtheorem{lemme}[theo]{Lemma}
\newtheorem{remark}[theo]{Remark}
\def \dem{\noindent{\sc\small Proof: }}
\def \findem {\hfill{\hbox {\vrule\vbox{\hrule width 6pt\vskip
6pt\hrule}\vrule}}}

\def\RM{{\mathbb{R}}}

\def\PM{{\mathbb{P}}}

\def\EM{{\mathbb{E}}}

\def \[{\text{\textlbrackdbl}}
\def \]{\text{\textrbrackdbl}}

\def \nrond{\mathcal{N}}

\def \frond{\mathcal{F}}
\def \hrond{\mathcal{H}}
\def \rrond{\mathcal{R}}

\def \urond{\mathcal{U}}
\def \prond{\mathcal{P}}

\begin{document}

\begin{frontmatter}



\title{{\LARGE {\bf Asymptotically efficient estimators for nonparametric heteroscedastic regression models}}\\}


\author{J.-Y. BRUA}
\ead{brua@math.u-strasbg.fr}
\address{IRMA, 7 rue Ren\'e Descartes F67084, Strasbourg Cedex, France}

\begin{abstract}
This paper concerns the estimation of a function at a point in
nonparametric heteroscedastic regression models with Gaussian
noise or noise having unknown distribution. In those cases an
asymptotically efficient kernel estimator is constructed for the
minimax absolute error risk.
\end{abstract}


\end{frontmatter}


\section{Introduction}

We consider the problem of estimating a regression function $S$ at
a given point $z_0\in]0;1[$ under observations
\begin{equation}\label{model}
y_k=S(x_k)+g(x_k,S)\xi_k,~~~~~~k\in\{1,\dots,n\}
\end{equation}
where the regressors $x_k=k/n$ are deterministic, $\xi_k$ are
independent identically distributed random variables which will
firstly be assumed Gaussian standard then having unknown density.
Notice that the variance of the noises $g^2$ is unknown and
depends on the unknown regression function $S$ and the regressors
$x_k$.

Heteroscedastic regression models with this type of scale
functionals have been encountered in consumer budget studies
utilizing observations on individuals with diverse incomes, in
analysis of investment behavior of firms of different sizes and
more recently in medical research. For example, \citet{goldfeld}
considered polynomial regression models such that
$y_k=\alpha+\beta x_k+u_k,~~\EM(u_k^2)=a+bx_k+cx_k^2,$ which is a
particular case of our model (\ref{model}) if we assume the
unknown regression function being like $S(x)=\alpha+\beta x$ and
$g^2(x,S)=(a-\frac{\alpha c}{\beta^2})+\left(b-2\frac{\alpha
c}{\beta}\right)x+\frac{c}{\beta^2}S^2(x)$. Other heteroscedastic
regression models are studied for instance in \citet{MR1422411},
\citet{galtchouk2} and \citet{MR2293406}.

The problem of H\"{o}lder regression estimation has been studied
by several authors. For a regression function belonging to a
quasi-H\"{o}lder class and estimated at a point with squared error
loss, \citet{MR606617} showed that the linear minimax estimator is
a kernel estimator. \citet{MR1105839} further found that this
estimator is within 17 percent of asymptotically minimax among all
procedures and obtained optimal kernels for H\"{o}lder classes.
For estimating the whole object or its \emph{k}th derivative with
sup-norm global loss and H\"{o}lder class, \citet{korostelev} and
\citet{MR1278880} proved that a kernel estimator is asymptotically
efficient.

This article deals with nonparametric estimation of a regression
function belonging to a H\"{o}lder ball. We work with the absolute
error loss and the corresponding risk. Our aim is to find an
efficient estimator, that is to say an estimator which achieves
the sharp asymptotic behavior of the minimax risk. To that purpose
we use the method developed by \citet{galtchouk} who introduce the
local weak H\"{o}lder classes to define the risk of an estimator.
So we use the classes $\urond_{z_0,\delta}$ which allows an
arbitrary large derivative but has a H\"{o}lder condition based on
a H\"{o}lder constant tending to zero (see (\ref{voisinage})),
then define the risk $\rrond_{z_0,\delta}(\tilde{S})$ of an
estimator $\tilde{S}$ of $S(z_0)$ and the minimax risk
$\displaystyle\inf_{\tilde{S}}\rrond_{z_0,\delta}(\tilde{S})$ (see
(\ref{risque})). In these conditions we prove that a kernel
estimator is asymptotically efficient, it means that the minimax
risk attains the sharp asymptotic constant.

This paper is organized as follows. In section 2 we describe the
problem in the case of Gaussian noise with all assumptions needed
and define all necessary mathematical objects. Our main results of
this problem are written in section 3. The case of unknown noise
is related in section 4. Theorems are proved in section 5 and
appendix A contains useful results for our proofs.

\section{Statement of the problem}

Consider model (\ref{model}) where $g:[0;1]\times
C^1([0;1],\RM)\longrightarrow\RM_+^*$ and $S$ are unknown
functions, $S$ belonging to the class
$$\hrond(\beta)=\bigcup_{M>0,K>0}\hrond (M,K,\beta),$$
where $\beta=1+\alpha$ is known, $\alpha\in ]0;1]$,
$\hrond(M,K,\beta)$ is the H\"{o}lder class defined as
$$\hrond (M,K,\beta)=\left\{S\in C^1([0;1],\RM)~:~\parallel S^{\prime}\parallel\leq M,~\sup_{x,y\in [0;1]}\frac{|S^{\prime}(y)-S^{\prime}(x)|}{|x-y|^{\alpha}}\leq K\right\},$$
with $\parallel f\parallel=\displaystyle\sup_{x\in [0;1]}|f(x)|$.
We suppose that the noises $(\xi_k)_{1\leq k\leq n}$ are
independent identically distributed $\nrond(0,1)$.

As mentioned in the introduction, we will work with a minimax risk
taken over the local weak H\"{o}lder class at the point $z_0$
defined, for $0<\delta<1$, as
\begin{equation}\label{voisinage}
\urond_{z_0,\delta}=\left\{S\in\hrond(\beta):\|
S^{\prime}\|\leq\delta^{-1};\forall
h>0,\left|\int_{-1}^{1}\big(S(z_0+hu)-S(z_0)\big)du\right|\leq\delta
h^{\beta}\right\}
\end{equation}
Notice that
\begin{equation}\label{condint}
\int_{-1}^{1}\big(S(z_0+hu)-S(z_0)\big)du=\int_{-1}^1\left(\int_{z_0}^{z_0+uh}(S^{\prime}(t)-S^{\prime}(z_0))dt\right)du,
\end{equation}
so we have for all $S\in\hrond(M,K,\beta)$
$$\left|\int_{-1}^{1}\big(S(z_0+hu)-S(z_0)\big)du\right|\leq
\frac{2K}{\beta(\beta+1)}h^{\beta}.$$ That is why the class
$\urond_{z_0,\delta}$ is called a weak H\"{o}lder class.\\
Furthermore (\ref{condint}) implies that
$\hrond(\delta^{-1},\delta,\beta)\subset\urond_{z_0,\delta}$ for
any $0<\delta<1$.

Let us give the assumptions needed. Firstly we suppose that
\begin{equation}\label{hyp1}
\lim_{n\rightarrow\infty}\sup_{S\in\urond_{z_0,\delta}}\left|\left(\frac{1}{q_n}\sum_{k=1}^{n}Q\big(\frac{x_k-z_0}{h}\big)g^2(x_k,S)\right)^{\frac{1}{2}}-g(z_0,S)\right|=0,
\end{equation}
with $$q_n=\sum_{k=1}^n Q\big(\frac{x_k-z_0}{h}\big),~~
Q=\mathbb{I}_{[-1;1]} \text{~and~} h=n^{-1/(2\beta+1)}.$$ Moreover
we assume that there exists $g_\star>0$ and $g^{\star}<\infty$
such that
\begin{equation}\label{gbornee}
g_{\star}\leq\inf_{0\leq x\leq 1}\inf_{S\in
C^1([0;1],\RM)}g(x,S)\leq\sup_{0\leq x\leq 1}\sup_{S\in
C^1([0;1],\RM)}g(x,S)\leq g^{\star}
\end{equation}
and that the function $g$ is differentiable in the Frechet sense
with respect to $S$ in $C^1([0;1],\RM)$ uniformly over
$x\in[0;1]$,
i.e. for any $S,S_0\in C^1([0;1],\RM)$\\
\begin{equation}\label{hyp2}
g(x,S)=g(x,S_0)+L_{x,S_0}(S-S_0)+\Gamma_{x,S_0}(S-S_0),
\end{equation}
where the linear operator $L_{x,S_0}$ is bounded on
$C^1([0;1],\RM)$ uniformly over $x\in[0;1]$, i.e. for any $S_0\in
C^1([0;1],\RM)$ there exists some positive constant $C_{S_0}$ such
that
\begin{equation}\label{hyp3}
\sup_{x\in[0;1]}\sup_{S\in C^1([0;1],\RM),\ \|S\|\neq
0}|L_{x,S_0}(S)|/\parallel S\parallel\leq C_{S_0}
\end{equation}
and the residual term $\Gamma_{x,S_0}(S)$ satisfies the property
\begin{equation}\label{hyp4}
\lim_{\parallel S\parallel\rightarrow
0}\sup_{x\in[0;1]}\Gamma_{x,S_0}(S)/\parallel S\parallel=0.
\end{equation}

\begin{remark}~\\
Note that hypothesis (\ref{hyp1}) is verified when for all
$\varepsilon>0$, there exists $\eta>0$ such that if
$|x-z_0|\leq\eta$, then
$\displaystyle\sup_{S\in C^1([0;1],\RM)}|g(x,S)-g(z_0,S)|\leq\varepsilon$.\\
In particular a function $g$ satisfies this property if it is
uniformly continuous with respect to both variables.
\end{remark}

\begin{remark}~\\
Let us give a general example of a function $g$ satisfying
hypothesis (\ref{hyp1})--(\ref{hyp4}) above. Let
$V:\RM\longrightarrow\RM_+$ and
$G:[0;1]\times\RM\longrightarrow\RM_+$ two differentiable
functions such that
$$\|V^{\prime}\|_{\infty}<\infty,~
G_{\star}=\inf_{x\in[0,1],\
y\in\RM}G(x,y)>0,~G_{\star}^{\prime}=\sup_{x\in[0;1],\
y\in\RM}\left|\frac{\partial G}{\partial y}(x,y)\right|<\infty.$$
Define
\begin{equation}\label{fonctiong}
g^2(x,S)=G(x,S(x))+\displaystyle\int_0^1V(S(t))dt.
\end{equation}
The derivative in the Frechet sense of $g$ is given by
$$L_{x,S}(f)=\frac{1}{2g(x,S)}\frac{\partial G}{\partial
y}(x,S(x))f(x)+\frac{1}{2g(x,S)}\int_0^1V^{\prime}(S(t))f(t)dt,$$
so we have
$$\sup_{x\in [0;1]}\sup_{S\in C^1([0;1],\RM),\ \|S\|\neq
0}\frac{|L_{x,S}(f)|}{\|f\|_{\infty}}\leq\frac{G_{\star}^{\prime}+\|V\|_{\infty}}{2\sqrt{G_{\star}}}.$$
Writing Taylor's expansion of functions $y\mapsto G(x,y)$ at the
point $(x,S(x))$ and $V$ at the point $S(t)$ to the first order:
\begin{eqnarray*}
G(x,S(x)+f(x)) & = & G(x,S(x))+\frac{\partial G}{\partial
y}(x,S(x))f(x)+f(x)\varepsilon_{x,S}(f(x)),\\
V(S(t)+f(t)) & = &
V(S(t))+V^{\prime}(S(t))f(t)+f(t)\tilde{\varepsilon}_{t,S}(f(t)),
\end{eqnarray*}
we can easily show that
$$\frac{|\Gamma_{x,S}(f)|}{\|f\|_{\infty}}\leq\frac{G_{\star}^{\prime}+\|V^{\prime}\|_{\infty}}{8G_{\star}^{3/2}}\left|g^2(x,S+f)-g^2(x,S)\right|$$
\begin{equation}\label{Gamma}
~~~~~~~~~~~~~~~~+\frac{1}{2\sqrt{G_{\star}}}\left(|\varepsilon_{x,S}(f(x))|+\int_0^1|\tilde{\varepsilon}_{t,S}(f(t))|dt\right).
\end{equation}
Now if we take $G(x,y)=\alpha_0+\alpha_1x+\alpha_2\sin^2y$ and
$V(y)=\alpha_3\sin^2y$ for all $(x,y)\in [0,1]\times\RM$, with
$\alpha_0>0$ and $\alpha_1,\alpha_2,\alpha_3\in\RM_+$, then the
function $g$ defined as (\ref{fonctiong}) is uniformly continuous,
bounded by $\sqrt{\alpha_0}$ and
$\sqrt{\alpha_0+\alpha_1+\alpha_2+\alpha_3}$. Moreover by writing
explicitly the functions $\varepsilon_{x,S}$ and
$\tilde{\varepsilon}_{x,S}$ for this case, we can prove thanks to
(\ref{Gamma}) that $g$ satisfies hypothesis (\ref{hyp4}). So we
have exhibited an example of function $g$ which satisfies all
assumptions needed.
\end{remark}

For any estimator $\tilde{S}_n(z_0)$ of $S(z_0)$ we define the
following risk
\begin{equation}\label{risque}
\rrond_{z_0,\delta}(\tilde{S}_n)=\sup_{S\in\urond_{z_0,\delta}}\EM_S\varphi_n\frac{|\tilde{S}_n(z_0)-S(z_0)|}{g(z_0,S)},
\end{equation}
where $\EM_S$ is the expectation taken  with respect to the law
$\PM_S$ in $(\ref{model})$ and
$\varphi_n=n^{\frac{\beta}{2\beta+1}}$.\\
The aim is to attain the sharp constant with this rate
$\varphi_n$. It is only assumed that $\beta\in]1;2]$ because if
$\beta>2$ we should use a kernel $Q$ of order $[\beta]$ i.e. such
that $\int u^jQ(u)du=0$ for $j=1,2,\dots,[\beta]$ and $\int
Q(u)du<\infty$, where $[a]$ denotes the integer part of the number
$a$.

\section{Lower and upper bounds}

In this section we give the lower bound for the minimax risk and
show that the kernel estimator $\hat{S}_n(z_0)$, defined by
\begin{equation}\label{estimateur}
\hat{S}_n(z_0)=\frac{1}{q_n}\sum_{k=1}^n
Q\big(\frac{x_k-z_0}{h}\big)y_k,
\end{equation}
is asymptotically efficient as we give the upper bound for its
risk.

\begin{theo}\label{thm3}
For any $\delta\in]0;1[$,
$$\liminf_{n\rightarrow\infty}\inf_{\tilde{S}}\rrond_{z_0,\delta}(\tilde{S})\geq\frac{\EM|\xi|}{\sqrt{2}},~~~~~~~\xi\sim\nrond(0,1),$$
where the infimum is taken over all estimators $\tilde{S}$ of
$S(z_0)$.
\end{theo}

\begin{theo}\label{thm4}
For the estimator $\hat{S}_n(z_0)$ from (\ref{estimateur}), the
following inequality holds:
$$\limsup_{\delta\rightarrow
0}\limsup_{n\rightarrow\infty}\rrond_{z_0,\delta}(\hat{S}_n(z_0))\leq\frac{\EM|\xi|}{\sqrt{2}},~~~~~~\xi\sim\nrond(0,1).$$
\end{theo}

\section{Case of unknown noise distribution}

In this section we suppose that the $(\xi_k)$ in model
(\ref{model}) are independent identically distributed with an
unknown density $p$ belonging to
$$\prond_{\varepsilon,L}=\left\{p: \int_{-\infty}^{+\infty}xp(x)dx=0, \int_{-\infty}^{+\infty}x^2p(x)dx=1, \int_{-\infty}^{+\infty}|x|^{2+\varepsilon}p(x)dx\leq L\right\},$$
with $\varepsilon>0$ and $L>0$ sufficiently large to have the
density of the standard Gaussian random variable in
$\prond_{\varepsilon,L}$.\\
We define the risk corresponding to this case as
$$\tilde{\rrond}_{z_0,\delta}(\tilde{S}_n)=\sup_{p\in\prond_{\epsilon,L}}\sup_{S\in\urond_{z_0,\delta}}\EM_S\varphi_n\frac{|\tilde{S}_n(z_0)-S(z_0)|}{g(z_0,S)}.$$
In the following theorems we give the sharp lower bound for the
minimax risk over all estimators and establish the upper bound for
the minimax risk for the kernel estimator $\hat{S}_n(z_0)$ of
$S(z_0)$ defined in (\ref{estimateur}).

\begin{theo}\label{thm1}
For any $\delta\in]0;1[$, one has:
$$\liminf_{n\rightarrow\infty}\inf_{\tilde{S}}\tilde{\rrond}_{z_0,\delta}(\tilde{S})\geq\frac{\EM|\eta|}{\sqrt{2}},~~~~~~\eta\sim\nrond(0,1),$$
where the infimum is taken over all estimators $\tilde{S}$ of
$S(z_0)$.
\end{theo}

\begin{theo}\label{thm2} The kernel estimator (\ref{estimateur}) is
asymptotically efficient. Indeed it satisfies the inequality:
$$\limsup_{\delta\rightarrow
0}\limsup_{n\rightarrow\infty}\tilde{\rrond}_{z_0,\delta}(\hat{S}_n(z_0))\leq\frac{\EM|\eta|}{\sqrt{2}},~~~~~~\eta\sim\nrond(0,1).$$
\end{theo}

\section{Proof of the theorems}
\subsection{Proof of theorem \ref{thm3}}
For all $\nu\in\left]0;\frac{1}{4}\right[$, denote
$S_{\nu}(x)=\varphi_n^{-1}V_{\nu}\left(\displaystyle\frac{x-z_0}{h}\right),$
where the function $V_{\nu}$ is defined by:
$$V_{\nu}(x)=\frac{1}{\nu}\int_{-\infty}^{+\infty}\tilde{Q}_{\nu}(u)l\left(\frac{u-x}{\nu}\right)du~,~
\tilde{Q}_{\nu}(u)=\mathbb{I}_{\{|u|\leq
1-2\nu\}}+2\mathbb{I}_{\{1-2\nu\leq|u|\leq 1-\nu\}},$$ and $l$ is
a non-negative function, infinitely differentiable on $\RM$, such
that for $|z|\geq 1$, $l(z)=0$ and
$\displaystyle\int_{-1}^1l(z)dz=1$. One can easily see that for
any $0<\nu<\frac{1}{4}$, we have $V_{\nu}(0)=1$ and
$\displaystyle\int_{-1}^1V_{\nu}(x)dx=2$.

Let $\nu\in]0;\frac{1}{4}[$, $b>0$ and $\delta\in]0;1[$. Denote
$S_{\nu,u}(x)=\displaystyle\frac{u}{\varphi_n}V_{\nu}\left(\displaystyle\frac{x-z_0}{h}\right),$
where $x,u\in\RM$.\\
Thanks to lemma \ref{lem1}, if $|u|\leq b$ there exists an integer
$n_{\nu,b,\delta}>0$ such that $S_{\nu,u}\in\urond_{z_0,\delta}$
for all $n\geq n_{\nu,b,\delta}$. Therefore for $n\geq
n_{\nu,b,\delta}$, one has:
\begin{eqnarray*}
\rrond_{z_0,\delta}(\tilde{S}) & \geq & \sup_{|u|\leq
b}\frac{1}{g(z_0,S_{\nu,u})}\EM_{S_{\nu,u}}\varphi_n|\tilde{S}(z_0)-S_{\nu,u}(z_0)|\\
& \geq &
\frac{1}{2b}\int_{-b}^b\frac{1}{g(z_0,S_{\nu,u})}\EM_{S_{\nu,u}}v_a\left(\varphi_n(\tilde{S}(z_0)-S_{\nu,u}(z_0))\right)du:=I_n(a,b),
\end{eqnarray*}
where $v_a(x)=|x|\wedge a, a>0$.

Write $\mathbb{P}_{S_{\nu,u}}$ the law of
$(y_k^{(1)})_{k=1,\dots,n}$, where
$y_k^{(1)}=S_{\nu,u}(x_k)+g(x_k,S_{\nu,u})\xi_k$, and $\mathbb{P}$
the law of $(y_k^{(0)})_{k=1,\dots,n}$, where
$y_k^{(0)}=g(x_k,S_{\nu,u})\xi_k$. These two measures are
equivalent and the corresponding Radon-Nikodym derivative is at
the point $(y_1,\dots,y_n)$:
\begin{eqnarray*}
\rho_n(u) & = &
\frac{d\PM_{S_{\nu,u}}}{d\PM}(y_1,\dots,y_n)\\
& = &
\exp\left\{-\frac{1}{2}\sum_{k=1}^n\left(\left(\frac{y_k-S_{\nu,u}(x_k)}{g(x_k,S_{\nu,u})}\right)^2-\left(\frac{y_k}{g(x_k,S_{\nu,u})}\right)^2\right)\right\}\\
& = &
\exp\left(u\varsigma_n\eta_n-\frac{u^2}{2}\varsigma_n^2\right)
\end{eqnarray*}
where
$\varsigma_n^2=\displaystyle\frac{1}{\varphi_n^2}\displaystyle\sum_{k=1}^n\frac{V_{\nu}^2\left(\frac{x_k-z_0}{h}\right)}{g^2(x_k,S_{\nu,u})}$
and
$\eta_n=\displaystyle\frac{1}{\varsigma_n\varphi_n}\displaystyle\sum_{k=1}^n\frac{V_{\nu}\left(\frac{x_k-z_0}{h}\right)}{g^2(x_k,S_{\nu,u})}y_k$.\\
Under the law $\mathbb{P}$, $\eta_n$ is a standard Gaussian random
variable.

We prove in lemma \ref{lem2} that
\begin{equation}
\label{sigmanu}\varsigma_n^2\xrightarrow[n\rightarrow\infty]{}\int_{-1}^1\frac{V_{\nu}^2(z)}{g^2(z_0,0)}dz=:\sigma_{\nu}^2.
\end{equation}
So we rewrite
$\rho_n(u)=\exp\left(u\sigma_{\nu}\eta_n-\frac{u^2\sigma_{\nu}^2}{2}+r_n\right)$,
where $r_n$ converges in $\PM$-probability to zero.

Denoting
$\psi_{a,n}(\tilde{S},S_{\nu,u})=v_a(\varphi_n(\tilde{S}_n(z_0)-S_{\nu,u}(z_0)))$
and $\EM$ the expectation for the probability measure $\PM$, one
has
\begin{equation}\label{JetDelta}
I_n(a,b)\geq\frac{1}{2b}\int_{-b}^b\EM\mathbb{I}_{B_d}\frac{\psi_{a,n}(\tilde{S},S_{\nu,u})}{g(z_0,S_{\nu,u})}\varrho_n(u)du+\delta_n(a,b)=:J_n(a,b)+\delta_n(a,b),
\end{equation}
where
\begin{eqnarray*}
B_d & = & \{|\eta_n|\leq d\}\text{~~and~~}
d=\sigma_{\nu}(b-\sqrt{b}), b>1,\\
\varrho_n(u) & = &
\exp\left(u\sigma_{\nu}\eta_n-\frac{u^2\sigma_{\nu}^2}{2}\right),\\
\delta_n(a,b) & = &
\frac{1}{2b}\int_{-b}^b\EM\mathbb{I}_{B_d}\frac{\psi_{a,n}(\tilde{S},S_{\nu,u})}{g(z_0,S_{\nu,u})}\theta_n(u)du,\\
\theta_n(u) & = & \rho_n(u)-\varrho_n(u).
\end{eqnarray*}

Note that
$\rho_n(u)\xrightarrow[n\rightarrow\infty]{\mathcal{L}}\rho_{\infty}(u)=\exp\left(u\sigma_{\nu}\eta-\frac{u^2\sigma_{\nu}^2}{2}\right)$.
We can easily show that $\EM\rho_{\infty}(u)=1$ and we have also
$\EM\rho_n(u)=1$ because $\rho_n(u)$ is a density. Hence, using
theorem 3.6 from \citet{billingsley}, $\{\rho_n(u), n\geq 1\}$ is
uniformly integrable. And since $\varrho_n(u)$ is bounded on
$B_d$, we obtain the uniform integrability of
$\{\mathbb{I}_{B_d}\psi_{a,n}(\tilde{S},S_{\nu,u})\theta_n(u),
n\geq 1\}$.\\
Write
$\theta_n(u)=\exp\left(u\sigma_{\nu}\eta_n-\frac{u^2\sigma_{\nu}^2}{2}\right)(e^{r_n}-1)$
and notice that
$\exp\left(u\sigma_{\nu}\eta_n-\frac{u^2\sigma_{\nu}^2}{2}\right)$
is bounded on $B_d$ and that
$e^{r_n}-1\displaystyle\xrightarrow[n\rightarrow\infty]{\mathbb{P}}0$.
As a consequence one has
$$\frac{\mathbb{I}_{B_d}\psi_{a,n}(\tilde{S},S_{\nu,u})}{g(z_0,S_{\nu,u})}\theta_n(u)\xrightarrow[n\rightarrow\infty]{\PM}0.$$
It follows that
$\frac{\mathbb{I}_{B_d}\psi_{a,n}(\tilde{S},S_{\nu,u})}{g(z_0,S_{\nu,u})}\theta_n(u)\xrightarrow[n\rightarrow\infty]{\mathbb{L}^1}0$
and $\EM\frac{\mathbb{I}_{B_d}\psi_{a,n}(\tilde{S},S_{\nu,u})}{g(z_0,S_{\nu,u})}\theta_n(u)\xrightarrow[n\rightarrow\infty]{}0$.\\
Finally bounded convergence yields
$\delta_n(a,b)\xrightarrow[n\rightarrow\infty]{}0$ in
(\ref{JetDelta}).

Now we are interested in the term $J_n(a,b)$ in
(\ref{JetDelta}).\\
First rewrite
$\varrho_n(u)=\zeta_ne^{-\sigma_{\nu}^2(u-\tilde{\eta}_n)^2/2}$
with $\zeta_n=e^{\eta_n^2/2}$ and
$\tilde{\eta}_n=\displaystyle\frac{\eta_n}{\sigma_{\nu}}$. Then if
$\xi\sim\nrond(0,1)$ denote
$\tilde{\xi}=\displaystyle\frac{\xi}{\sigma_{\nu}}$,
$\zeta=e^{\xi^2/2}$, $\tilde{B}_d=\{|\xi|\leq d\}$ and
$\tilde{\EM}$ the expectation for the probability law of $\xi$.
With $t_n=\varphi_n\tilde{S}_n(z_0)$, we get
\begin{eqnarray*}
J_n(a,b) & = &
\frac{1}{2b}\int_{-b}^b\EM\mathbb{I}_{B_d}\zeta_n\frac{v_a(u-t_n)}{g(z_0,S_{\nu,u})}\exp\left(-\frac{\sigma_{\nu}^2}{2}(u-\tilde{\eta}_n)^2\right)du\\
& = &
\frac{1}{2b}\int_{-b}^b\tilde{\EM}\mathbb{I}_{\tilde{B}_d}\zeta\frac{v_a(u-t_n)}{g(z_0,S_{\nu,u})}\exp\left(-\frac{\sigma_{\nu}^2}{2}(u-\tilde{\xi})^2\right)du\\
& = &
\tilde{\EM}\mathbb{I}_{\tilde{B}_d}\zeta\frac{1}{2b}\int_{-b}^b\frac{v_a(u-t_n)}{g(z_0,S_{\nu,u})}\exp\left(-\frac{\sigma_{\nu}^2}{2}(u-\tilde{\xi})^2\right)du.
\end{eqnarray*}
We have the following limit
\begin{equation}\label{limJn}
\tilde{\EM}\mathbb{I}_{\tilde{B}_d}\zeta\frac{1}{2b}\int_{-b}^bv_a(u-t_n)\exp\left(-\frac{\sigma_{\nu}^2}{2}(u-\tilde{\xi})^2\right)\left(\frac{1}{g(z_0,S_{\nu,u})}-\frac{1}{g(z_0,0)}\right)du\xrightarrow[n\rightarrow\infty]{}0.
\end{equation}
Indeed, using hypothesis (\ref{hyp2}) and (\ref{hyp3}) one obtains
\begin{eqnarray*}
& &
\left|\tilde{\EM}\mathbb{I}_{\tilde{B}_d}\zeta\frac{1}{2b}\int_{-b}^bv_a(u-t_n)\exp\left(-\frac{\sigma_{\nu}^2}{2}(u-\tilde{\xi})^2\right)\left(\frac{1}{g(z_0,S_{\nu,u})}-\frac{1}{g(z_0,0)}\right)du\right|\\
& \leq &
\tilde{\EM}\mathbb{I}_{\tilde{B}_d}\zeta\frac{1}{2b}\int_{-b}^bv_a(u-t_n)\exp\left(-\frac{\sigma_{\nu}^2}{2}(u-\tilde{\xi})^2\right)\left|\frac{\Gamma_{z_0,0}(S_{\nu,u})-L_{z_0,0}(S_{\nu,u})}{g_{\star}^2}\right|du\\
& \leq &
\tilde{\EM}\mathbb{I}_{\tilde{B}_d}\zeta\frac{1}{2b}\int_{-b}^ba\frac{C_0\parallel S_{\nu,u}\parallel+|\Gamma_{z_0,0}(S_{\nu,u})|}{g_{\star}^2}du.\\
\end{eqnarray*}
Since $\parallel S_{\nu,u}\parallel$ tends to zero as $n$ goes to
infinity, hypothesis (\ref{hyp4}) and (\ref{limJn}) allows then us
to say that
$$\liminf_{n\rightarrow\infty}J_n(a,b)=\liminf_{n\rightarrow\infty}\tilde{\EM}\mathbb{I}_{\tilde{B}_d}\zeta\frac{1}{2b}\int_{-b}^b\frac{v_a(u-t_n)}{g(z_0,0)}\exp\left(-\frac{\sigma_{\nu}^2}{2}(u-\tilde{\xi})^2\right)du.$$
But \begin{eqnarray*} & &
\tilde{\EM}\mathbb{I}_{\tilde{B}_d}\zeta\frac{1}{2b}\int_{-b}^b\frac{v_a(u-t_n)}{g(z_0,0)}\exp\left(-\frac{\sigma_{\nu}^2}{2}(u-\tilde{\xi})^2\right)du\\
& \geq &
\tilde{\EM}\mathbb{I}_{\tilde{B}_d}\zeta\frac{1}{2b}\int_{-\sqrt{b}}^{\sqrt{b}}\frac{v_a(t-t_n+\tilde{\xi})}{g(z_0,0)}\exp\left(-\frac{\sigma_{\nu}^2}{2}t^2\right)dt\\
& \geq &
\tilde{\EM}\mathbb{I}_{\tilde{B}_d}\zeta\frac{1}{2b}\int_{-\sqrt{b}}^{\sqrt{b}}\frac{v_a(t)}{g(z_0,0)}\exp\left(-\frac{\sigma_{\nu}^2}{2}t^2\right)dt,
\end{eqnarray*}
this last inequality holds thanks to Anderson's lemma \citep[see][Chapter II, Lemma 10.1 and Corollary 10.2]{ibragimov}.\\
Eventually using the fact that
$\tilde{\EM}\mathbb{I}_{\tilde{B}_d}\zeta=\displaystyle\frac{2\sigma_{\nu}(b-\sqrt{b})}{\sqrt{2\pi}}$
it follows that
$$\liminf_{a\rightarrow\infty}\liminf_{n\rightarrow\infty}J_n(a,b)\geq
\frac{\sigma_{\nu}}{\sqrt{2\pi}}\frac{b-\sqrt{b}}{b}\int_{-\sqrt{b}}^{\sqrt{b}}\frac{|t|}{g(z_0,0)}\exp\left(-\frac{\sigma_{\nu}^2}{2}t^2\right)dt.$$
We complete the proof limiting successively $b\rightarrow\infty$,
$\nu\rightarrow 0$ and utilizing
$\sigma_{\nu}^2\xrightarrow[\nu\rightarrow
0]{}\displaystyle\frac{2}{g^2(z_0,0)}$.\findem

\subsection{Proof of theorem \ref{thm4}}
We begin by rewriting the kernel estimator as $\hat{S}_n(z_0)-S(z_0)=B_n+\frac{1}{\sqrt{q_n}}\zeta_n$ with\\
\begin{eqnarray}
B_n & = & \frac{1}{q_n}\sum_{k=1}^{n}Q\big(\frac{x_k-z_0}{h}\big)(S(x_k)-S(z_0))\label{Bn}\\
\zeta_n & = &
\frac{1}{\sqrt{q_n}}\sum_{k=1}^{n}Q\big(\frac{x_k-z_0}{h}\big)g(x_k,S)\xi_k.\label{zetan}
\end{eqnarray}

First we take a look at the term
$\displaystyle\frac{\zeta_n}{\sqrt{q_n}}$. By (\ref{zetan}),
$\zeta_n$ is a Gaussian random variable
$\nrond\big(0,\sigma_n^2(S)\big)$ where
$\sigma_n^2(S)=\displaystyle\frac{1}{q_n}\sum_{k=1}^nQ\big(\frac{x_k-z_0}{h}\big)g^2(x_k,S)$.
We prove in lemma \ref{variance} that the variance $\sigma_n^2(S)$
satisfies
$\sigma_n^2(S)\xrightarrow[n\rightarrow\infty]{}g^2(z_0,S)$. If
$\xi\sim\nrond(0,1)$, one has
\begin{eqnarray*}
\sup_{S\in\urond_{z_0,\delta}}\frac{1}{g(z_0,S)}\EM_S\big|\frac{\varphi_n}{\sqrt{q_n}}\zeta_n\big|
& = &
\frac{\varphi_n}{\sqrt{q_n}}\EM|\xi|\sup_{S\in\urond_{z_0,\delta}}\frac{\sigma_n(S)}{g(z_0,S)}\\
& \leq & \frac{\varphi_n}{\sqrt{q_n}}\frac{\EM|\xi|}{g_{\star}}\left(\sup_{S\in\urond_{z_0,\delta}}\left|\sigma_n(S)-g(z_0,S)\right|+g_{\star}\right).\\
\end{eqnarray*}
According to hypothesis $(\ref{hyp1})$ and since
$\displaystyle\frac{q_n}{\varphi_n^2}=\frac{q_n}{nh}\xrightarrow[n\rightarrow\infty]{}2$,
we obtain
\begin{equation}\label{zeta}
\limsup_{n\rightarrow\infty}\sup_{S\in\urond{z_0,\delta}}\EM_S\frac{\varphi_n}{g(z_0,S)}\frac{|\zeta_n|}{\sqrt{q_n}}\leq\frac{\EM|\xi|}{\sqrt{2}}.
\end{equation}

Now denote $u_k=\displaystyle\frac{x_k-z_0}{h}$, $\Delta
u_k=\displaystyle\frac{1}{nh}$ and rewrite (\ref{Bn}) as
\begin{eqnarray}
B_n & = & \frac{\varphi_n^2}{q_n}\sum_{k=1}^nQ(u_k)\big(S(z_0+hu_k)-S(z_0)\big)\Delta u_k\\
& = &
\frac{\varphi_n^2}{q_n}\int_{-1}^1\big(S(z_0+hu)-S(z_0)\big)du+\frac{\varphi_n^2}{q_n}R_n\label{autre}
\end{eqnarray}
with
\begin{eqnarray*}
R_n & = & \sum_{k=1}^nQ(u_k)\big(S(z_0+hu_k)-S(z_0)\big)\Delta u_k-\int_{-1}^1\big(S(z_0+hu)-S(z_0)\big)du\\
& = & \sum_{k=k_{*}}^{k^{*}}\int_{u_{k-1}}^{u_k}\big(S(z_0+hu_k)-S(z_0+hu)\big)du\\
& - &
\int_{u_{k^{*}}}^1\big(S(z_0+hu)-S(z_0)\big)du+\int_{u_{k_{*}-1}}^{-1}\big(S(z_0+hu)-S(z_0)\big)du,
\end{eqnarray*}
where $k^{*}=[n(z_0+h)]$ et $k_{*}=[n(z_0-h)]+1$.

We can bound $R_n$ as follows:
\begin{eqnarray*}
|R_n| & \leq & \sum_{k=k_{*}}^{k^{*}}\int_{u_{k-1}}^{u_k}h(u_k-u)\delta^{-1}du+\int_{u_{k^{*}}}^1h\delta^{-1}udu +\int_{u_{k_{*}-1}}^{-1}h\delta^{-1}|u|du\\
& \leq &
h\delta^{-1}\Big(\sum_{k=k_{*}}^{k^{*}}(u_k-u_{k-1})\frac{1}{nh}+(1-u_{k^{*}})+2(-1-u_{k_{*}-1})\Big)\leq
\frac{6\delta^{-1}}{n}.
\end{eqnarray*}
Hence
\begin{equation}\label{Rn}
\limsup_{n\rightarrow\infty}\sup_{S\in\urond_{z_0,\delta}}\EM_S\varphi_n\Big|\frac{\varphi_n^2}{q_n}R_n\Big|=0.
\end{equation}

With regard to the term
$\displaystyle\frac{\varphi_n^2}{q_n}\displaystyle\int_{-1}^1\big(S(z_0+hu)-S(z_0)\big)du$
in (\ref{autre}) one has
$$\Big|\frac{\varphi_n^2}{q_n}\int_{-1}^1\big(S(z_0+hu)-S(z_0)\big)du\Big|\leq\displaystyle\frac{\varphi_n^2}{q_n}\delta n^{\frac{-\beta}{2\beta+1}}=\delta\frac{\varphi_n}{q_n}.$$
Then using the definition of $\urond_{z_0,\delta}$ we get
\begin{equation}\label{delta}
\limsup_{n\rightarrow\infty}\sup_{S\in\urond_{z_0,\delta}}\EM_S\varphi_n\Big|\frac{\varphi_n^2}{q_n}\int_{-1}^1\big(S(z_0+hu)-S(z_0)\big)du\Big|
\leq\frac{\delta}{2}.
\end{equation}
Finally (\ref{zeta}), (\ref{Rn}) and limiting $\delta\rightarrow
0$ in (\ref{delta}) yield
$$\limsup_{\delta\rightarrow 0}\limsup_{n\rightarrow\infty}\rrond_{z_0,\delta}(\hat{S}_n(z_0))\leq\frac{\EM|\xi|}{\sqrt{2}}.$$\findem

\subsection{Proof of theorem \ref{thm1}}

This is a consequence of the theorem \ref{thm3} which gives the
sharp lower bound in the case of Gaussian errors having
expectation zero and unknown variance which depends on the design
point and the regression function. The corresponding risk
$\rrond_{z_0,\delta}$ is less than the risk
$\tilde{\rrond}_{z_0,\delta}$ because the density of the standard
Gaussian random variable belongs to $\prond_{\varepsilon,L}$. The
inequality in theorem \ref{thm1} is then proved.\findem

\subsection{Proof of theorem \ref{thm2}}

Writing $\hat{S}_n(z_0)-S(z_0)=B_n+\zeta_n/\sqrt{q_n}$, with $B_n$
and $\zeta_n$ defined by (\ref{Bn}) and (\ref{zetan}), we remark
that $B_n$ does not depend on the distributions of the random
variables $\xi_k$. That is the reason why (\ref{Rn}) and
(\ref{delta}) remain available and provide for any
$\delta\in]0;1[$:
$$\limsup_{n\rightarrow\infty}\sup_{S\in\urond_{z_0,\delta}}\varphi_n|B_n|\leq\delta/2.$$
Hence it suffices to prove that
\begin{equation}\label{esper}
\lim_{n\rightarrow\infty}\sup_{p\in\prond_{\epsilon,L}}\sup_{S\in\urond_{z_0,\delta}}\left|\frac{\EM_S|\zeta_n|}{g(z_0,S)}-\EM|\eta|\right|=0,
\end{equation}
with $\eta\sim\nrond(0,1)$.

Denote
$\tilde{\zeta_n}=\zeta_n/g(z_0,S)=\displaystyle\sum_{k=1}^nu_k$,
where
$u_k=\displaystyle\frac{1}{\sqrt{q_n}}Q\left(\frac{x_k-z_0}{h}\right)\frac{g(x_k,S)}{g(z_0,S)}\xi_k$,
and rewrite
$\displaystyle\frac{g(x_k,S)}{g(z_0,S)}\xi_k=\xi_k^{\prime}+\xi_k^{\prime\prime}$,
where
\begin{eqnarray*}
\xi_k^{\prime} & = &
\frac{g(x_k,S)}{g(z_0,S)}\xi_k\mathbb{I}_{|\xi_k|\leq
q_n^{1/4}}-\frac{g(x_k,S)}{g(z_0,S)}\EM\left(\xi_1\mathbb{I}_{|\xi_1|\leq q_n^{1/4}}\right),\\
\xi_k^{\prime\prime} & = &
\frac{g(x_k,S)}{g(z_0,S)}\xi_k\mathbb{I}_{|\xi_k|>
q_n^{1/4}}-\frac{g(x_k,S)}{g(z_0,S)}\EM\left(\xi_1\mathbb{I}_{|\xi_1|>
q_n^{1/4}}\right).
\end{eqnarray*}
Let
$u_k^{\prime}=\displaystyle\frac{1}{\sqrt{q_n}}Q\left(\displaystyle\frac{x_k-z_0}{h}\right)\xi_k^{\prime}$
and
$u_k^{\prime\prime}=\displaystyle\frac{1}{\sqrt{q_n}}Q\left(\displaystyle\frac{x_k-z_0}{h}\right)\xi_k^{\prime\prime}$,
then one gets
$\tilde{\zeta}_n=\tilde{\zeta}_n^{\prime}+\tilde{\zeta}_n^{\prime\prime}=\displaystyle\sum_{k=1}^nu_k^{\prime}+\displaystyle\sum_{k=1}^nu_k^{\prime\prime}$.
Moreover, $(u_k^{\prime})_{k\geq 1}$ is a martingale difference
and for all $k\geq 2$, we have $|u_k^{\prime}|\leq
2\frac{g^{\star}}{g_{\star}}q_n^{-1/4}$ and
$$\EM_S\left((u_{k}^{\prime})^2|\frond_{k-1}\right)=
\frac{1}{q_n}Q\left(\frac{x_k-z_0}{h}\right)\frac{g^2(x_k,S)}{g^2(z_0,S)}Var\left(\xi_{1}\mathbb{I}_{|\xi_1|\leq
q_n^{1/4}}\right).$$
Write
$$\sum_{i=1}^{n}\EM_S\left((u_i^{\prime
})^2|\frond_{i-1}\right)=\frac{Var\left(\xi_{1}\mathbb{I}_{|\xi_1|\leq
q_n^{1/4}}\right)}{q_n}\sum_{i=1}^{n}Q\left(\frac{x_i-z_0}{h}\right)\frac{g^2(x_i,S)}{g^2(z_0,S)}=\frac{G_n(S)}{q_n}a_n,$$
where
$G_n(S)=\displaystyle\sum_{i=1}^{n}Q\left(\frac{x_i-z_0}{h}\right)\frac{g^2(x_i,S)}{g^2(z_0,S)}$
and $a_n=Var\left(\xi_{1}\mathbb{I}_{|\xi_1|\leq
q_n^{1/4}}\right)$.\\
Denoting $r_n(S)=\displaystyle\frac{G_n(S)}{q_n}a_n$ and
$\tau_n=\inf\left\{k :
\displaystyle\sum_{i=1}^k\EM_S\left(u_i^{\prime
2}|\frond_{i-1}\right)\geq r_n(S)\right\}$, we obtain
$\tau_n=\inf\left\{k :
\displaystyle\sum_{i=1}^kQ\left(\frac{x_i-z_0}{h}\right)\geq
q_n\right\}$ and
$\tilde{\zeta}_n^{\prime}=\displaystyle\sum_{k=1}^{\tau_n}u_k^{\prime}$.

Let us show that $a_n$ and further $r_n(S)$ tend to $1$ uniformly
in $p\in\prond_{\epsilon,L}$ and in $S\in\urond_{z_0,\delta}$.
Firstly we have:
\begin{eqnarray*}
|a_n-1| & = & |\EM\left(\xi_1^2\mathbb{I}_{|\xi_1|\leq
q_n^{1/4}}\right)-\EM\left(\xi_1\mathbb{I}_{|\xi_1|\leq
q_n^{1/4}}\right)^2-1|\\
& \leq &
\left|\int_{-q_n^{1/4}}^{q_n^{1/4}}x^2p(x)dx-1\right|+\left|\int_{-q_n^{1/4}}^{q_n^{1/4}}xp(x)dx\right|^2.
\end{eqnarray*}
The Cauchy-Schwarz inequality brings us:
$$\left|\int_{-q_n^{1/4}}^{q_n^{1/4}}xp(x)dx\right|^2\leq
\left(\int_{-\infty}^{+\infty}x^2p(x)\mathbb{I}_{|x|>q_n^{1/4}}dx\right)\left(\int_{-\infty}^{+\infty}p(x)dx\right)
\leq K_p(q_n^{1/4}).$$ Nevertheless by the definition of the set
$\prond_{\epsilon,L}$, we get
\begin{equation}\label{Kp}
\displaystyle\sup_{p\in\prond_{\epsilon,L}}K_p(a):=\displaystyle\sup_{p\in\prond_{\epsilon,L}}\int_{-\infty}^{+\infty}x^2\mathbb{I}_{|x|>a}p(x)dx\displaystyle\xrightarrow[a\rightarrow\infty]{}0.
\end{equation}
From here it follows that
$$\sup_{p\in\prond_{\epsilon,L}}\sup_{S\in\urond_{z_0,\delta}}|a_n-1|\leq
2\sup_{p\in\prond_{\epsilon,L}}K_p(q_n^{1/4}),$$ so the left term
goes to zero as $n$ goes to infinity.\\
Using assumption (\ref{hyp1}) and the inequality
$$|r_n(S)-1|\leq\left|\frac{G_n(S)}{q_n}-1\right|+\frac{G_n(S)}{q_n}|a_n-1|$$
we get the convergence of $r_n(S)$ to $1$ uniformly in $p$ and in
$S$.

Applying lemma \ref{freedman} shows on the one hand the
convergence in distribution of $\zeta_n^{\prime}$ to $\nrond(0,1)$
uniformly in $p\in\prond_{\epsilon,L}$ and in
$S\in\urond_{z_0,\delta}$ because the function $\rho$ in lemma
\ref{freedman} does not depend on the law of the martingale
difference. In fact, if $\Phi$ denotes the standard Gaussian
distribution function, one has
\begin{eqnarray*}
& &\left|\PM\left(\sum_{k=1}^{\tau_n}u_k^{\prime}\leq
x\right)-\Phi(x)\right|\\ & \leq &
\left|\PM\left(\sum_{k=1}^{\tau_n}u_k^{\prime}\leq
x\right)-\Phi(x/\sqrt{r_n(S)})\right| +
\left|\Phi(x)-\Phi(x/\sqrt{r_n(S)})\right|.
\end{eqnarray*}
The second term of the right member of this inequality tends
toward zero uniformly in $p$, in $S$and in $x$ because
$r_n(S)\rightarrow 1$ uniformly in $p$ and in $S$ and because
$\Phi$ is uniformly continuous on $\RM$.

On the other hand one has $\EM|\tilde{\zeta}_n^{\prime\prime
}|\rightarrow 0$ uniformly in $p$ and in $S$. Indeed one have
immediately $\EM(\tilde{\zeta}_n^{\prime\prime
2})=\displaystyle\frac{G_n(S)}{q_n}K_p(q_n^{1/4}).$ Then
(\ref{Kp}) and the Cauchy-Schwarz inequality yield
$$\sup_{p\in\prond_{\epsilon,L}}\sup_{S\in\urond_{z_0,\delta}}\EM_S|\tilde{\zeta}_n^{\prime\prime
}|\rightarrow 0.$$ Using Markov's inequality, we show that
$(\tilde{\zeta}_n^{\prime\prime})$ tends to $0$ in probability
uniformly in $p$ and in $S$.

As a consequence
$\tilde{\zeta}_n=\tilde{\zeta}_n^{\prime}+\tilde{\zeta}_n^{\prime\prime}$
converges in distribution to $\eta\sim\nrond(0,1)$ uniformly in
$p$ and in $S$. This immediately implies (\ref{esper}).\findem

\appendix
\section{Appendix}

\begin{lemme}\label{lem1}Fix $\nu\in]0;\frac{1}{4}[$ and $\delta\in]0;1[$. Then there exists an integer
$n_{\nu,\delta}>0$ such that $S_{\nu}\in\urond_{z_0,\delta}$ for
all $n\geq n_{\nu,\delta}$.
\end{lemme}
\dem First remark that
$\displaystyle\int_{-1}^1\left(S_{\nu}(z_0+uh)-S_{\nu}(z_0)\right)du=0$.
Moreover one has
$$|S_{\nu}^{\prime}(x)|=\frac{1}{\varphi_nh}\left|V_{\nu}^{\prime}\left(\frac{x-z_0}{h}\right)\right|\leq\frac{2\|l^{\prime}\|_{\infty}}{\nu^2}n^{\frac{-\beta+1}{2\beta+1}}$$
For any fixed $\delta$ in $]0;1[$, if we choose $n\geq 1$ such
that
$$n^{\frac{-\beta+1}{2\beta+1}}\frac{2\|l^{\prime}\|_{\infty}}{\nu^2}\leq\delta^{-1}\text{~~i.e.~~}n\geq\left(\frac{2\|l^{\prime}\|_{\infty}\delta}{\nu^2}\right)^{\frac{2\beta+1}{\beta-1}},$$
then $S_{\nu}\in\urond_{z_0,\delta}$.\\
Therefore we have the desired result.\findem

\begin{lemme}\label{lem2}
We have the following limit:
$$\varsigma_n^2\xrightarrow[n\rightarrow\infty]{}\int_{-1}^1\frac{V_{\nu}^2(z)}{g^2(z_0,0)}dz.$$
\end{lemme}
\dem For sufficiently large $n$ we have $$\varsigma_n^2=
\frac{1}{nh}\sum_{k=1}^n\frac{V_{\nu}^2\left(\frac{x_k-z_0}{h}\right)}{g^2(x_k,S_{\nu,u})}=\frac{1}{h}\int_{z_0-h}^{z_0+h}\frac{V_{\nu}^2\left(\frac{x-z_0}{h}\right)}{g^2(x,S_{\nu,u})}\mu_n(dx)=\int_0^1\frac{V_{\nu}^2\left(\frac{x-z_0}{h}\right)}{g^2(x,S_{\nu,u})}\nu_n(dx)$$
with $\mu_n=\frac{1}{n}\displaystyle\sum_{k=1}^n\delta_{k/n}=$ and
$\nu_n=\frac{\mathbb{I}_{[z_0-h;z_0+h]}}{h}\mu_n$.\\
Using hypothesis $(\ref{hyp2})$ and $(\ref{hyp3})$ to the function
$g$, we can write for all $x\in[0;1]$
\begin{eqnarray*}
\left|\frac{1}{g^2(x,S_{\nu,u})}-\frac{1}{g^2(x,0)}\right|& \leq &
\frac{1}{g_{\star}^4}\left|2g(x,0)L_{x,0}(S_{\nu,u})+L_{x,0}^2(S_{\nu,u})+\Gamma_{x,0}^2(S_{\nu,u})\right.\\
& + & \left. 2g(x,0)\Gamma_{x,0}(S_{\nu,u})+2L_{x,0}(S_{\nu,u})\Gamma_{x,0}(S_{\nu,u})\right|\\
& \leq & \frac{1}{g_{\star}^4}\left(2g^{\star}C_0\|
S_{\nu,u}\|+C_0^2\|
S_{\nu,u}\|^2+|\Gamma_{x,0}(S_{\nu,u})|^2\right.\\
& + & \left.
2g^{\star}|\Gamma_{x,0}(S_{\nu,u})|+2C_0\|S_{\nu,u}\||\Gamma_{x,0}(S_{\nu,u})|\right).
\end{eqnarray*}
Hence
\begin{eqnarray*}
& &\left|\int_0^1\left(\frac{1}{g^2(x,S_{\nu,u})}-\frac{1}{g^2(x,0)}\right)\nu_n(dx)\right|\\
& \leq &
\frac{\|S_{\nu,u}\|}{g_{\star}^4}\int_0^1\nu_n(dx)\left(2g^{\star}C_0+C_0^2\|S_{\nu,u}\|+\left(\sup_{x\in[0;1]}\frac{|\Gamma_{x,0}(S_{\nu,u})|}{\|
S_{\nu,u}\|}\right)^2\|
S_{\nu,u}\|\right.\\
& + & \left.
2g^{\star}\left(\sup_{x\in[0;1]}\frac{|\Gamma_{x,0}(S_{\nu,u})|}{\|
S_{\nu,u}\|}\right)+2C_0\left(\sup_{x\in[0;1]}\frac{|\Gamma_{x,0}(S_{\nu,u})|}{\|
S_{\nu,u}\|}\right)\| S_{\nu,u}\|\right).
\end{eqnarray*}
As $(\nu_n)$ weakly tends to $2\delta_{z_0}$ when
$n\rightarrow\infty$, one has
$$\lim_{n\rightarrow\infty}\displaystyle\int_0^1\nu_n(dx)=2\text{~~et~~}\lim_{n\rightarrow\infty}\int_0^1\left(\frac{1}{g^2(x,0)}-\frac{1}{g^2(z_0,0)}\right)\nu_n(dx)=0.$$
Then taking into account hypothesis $(\ref{hyp4})$ and because
$\parallel S_{\nu,u}\parallel$ tends to $0$ as
$n\rightarrow\infty$, we obtain on the one hand
\begin{eqnarray*}
\int_0^1\left(\frac{1}{g^2(x,S_{\nu,u})}-\frac{1}{g^2(z_0,0)}\right)\nu_n(dx)\xrightarrow[n\rightarrow\infty]{}0.
\end{eqnarray*}
On the other hand
$$\int_0^1\frac{V_{\nu}^2\left(\frac{x-z_0}{h}\right)}{g^2(z_0,0)}\nu_n(dx)
=\frac{1}{nh}\sum_{k=1}^n\frac{V_{\nu}^2\left(\frac{x_k-z_0}{h}\right)}{g^2(z_0,0)}
\xrightarrow[n\rightarrow\infty]{}\int_{-1}^1\frac{V_{\nu}^2(y)}{g^2(z_0,0)}dy.$$
Now, if $V_{\nu}^{\star}$ denotes the maximum of $V_{\nu}^2$ on
$\RM$, one has
\begin{eqnarray*}
\left|\varsigma_n^2-\int_{-1}^1\frac{V_{\nu}^2(z)}{g^2(z_0,0)}dz\right|
& \leq &
\left|\int_0^1\frac{V_{\nu}^2\left(\frac{x-z_0}{h}\right)}{g^2(x,S_{\nu,u})}\nu_n(dx)-\int_0^1\frac{V_{\nu}^2\left(\frac{x-z_0}{h}\right)}{g^2(z_0,0)}\nu_n(dx)\right|\\
& + & \left|\int_0^1\frac{V_{\nu}^2\left(\frac{x-z_0}{h}\right)}{g^2(z_0,0)}\nu_n(dx)-\int_{-1}^1\frac{V_{\nu}^2(z)}{g^2(z_0,0)}dz\right|\\
& \leq &
V_{\nu}^{\star}\int_0^1\left|\frac{1}{g^2(x,S_{\nu,u})}-\frac{1}{g^2(z_0,0)}\right|\nu_n(dx)\\
& + &
\left|\int_0^1\frac{V_{\nu}^2\left(\frac{x-z_0}{h}\right)}{g^2(z_0,0)}\nu_n(dx)-\int_{-1}^1\frac{V_{\nu}^2(z)}{g^2(z_0,0)}dz\right|.
\end{eqnarray*}
Let $n$ goes to $\infty$ and then we have completed the proof of
lemma \ref{lem2}.\findem

\begin{lemme}\label{variance}
The variance $\sigma_n^2(S)$ of $\zeta_n$ satisfies
$$\sigma_n^2(S)\xrightarrow[n\rightarrow\infty]{}g^2(z_0,S).$$
\end{lemme}
\dem One has
$$\sum_{k=1}^nQ\big(\frac{x_k-z_0}{h}\big)g^2(x_k,S)=n\int_{z_0-h}^{z_0+h}g^2(x,S)\mu_n(dx)$$
with the measure $\mu_n=\frac{1}{n}\sum_{k=1}^n\delta_{k/n}$.\\
We know that $(\mu_n)_{n\geq 1}$ weakly tends to the uniform measure on $[0;1]$.\\
Moreover for sufficiently large $n$,
$$\frac{1}{nh}\sum_{k=1}^nQ\big(\frac{x_k-z_0}{h}\big)g^2(x_k,S)=\int_{z_0-h}^{z_0+h}g^2(x,S)\nu_n(dx)$$
with $\nu_n=\frac{\mu_n\mathbb{I}_{[z_0-h;z_0+h]}}{h}$.\\
Like this $(\nu_n)_{n\geq 1}$ weakly tends to $2\delta_{z_0}$, the Dirac measure at $z_0$, when $n\rightarrow\infty$.\\
Then we can conclude as we remember that
$\displaystyle\frac{q_n}{\varphi_n^2}\xrightarrow[n\rightarrow\infty]{}2$
and that $nh=\varphi_n^2$.\findem

\begin{lemme}\label{freedman} \citep[][pp. 90-91]{freedman}
Let $\delta\in]0;1[$ and $r>0$. Assume that $(u_k)_{k\geq 0}$ is a
martingale difference with respect to the filtration
$(\frond_k)_{k\geq 0}$ such that $|u_k|\leq\delta$ for all $k$ and
$\displaystyle\sum_{k=1}^{\infty}\EM(u_k^2|\frond_{k-1})\geq r$.\\
Define $\tau=\inf\left\{n :
\displaystyle\sum_{k=1}^n\EM(u_k^2|\frond_{k-1})\geq r\right\}$.\\
Then there exists a function $\rho:~]0;+\infty[\rightarrow [0;2]$
not depending on the distribution of the martingale difference,
such that $\displaystyle\lim_{x\rightarrow 0}\rho(x)=0$ and
$$\sup_{x\in\RM}\left|\PM\left(\sum_{k=1}^{\tau}u_k\leq
x\right)-\Phi(x/\sqrt{r})\right|\leq\rho(\delta/\sqrt{r}),$$ where
$\Phi$ is the standard Gaussian distribution function.
\end{lemme}




\end{document}